\documentclass[a4paper,12pt]{article}
\usepackage{amsfonts,amssymb,latexsym,amsmath}
\usepackage{multicol}
\usepackage[cp1251]{inputenc}
\usepackage[russian]{babel}
\usepackage{rotating}
\usepackage{graphicx}
\usepackage{tabularx}
\begin{document}
\newcounter{bnomer}
\newcounter{snomer}
\newcounter{diagram}
\setcounter{bnomer}{0} \setcounter{diagram}{0}
\renewcommand{\thesnomer}{\thebnomer.\arabic{snomer}}
\renewcommand{\thebnomer}{\arabic{bnomer}}

\newcommand{\sect}[1]{%
\setcounter{snomer}{0} \refstepcounter{bnomer}
\begin{center}\large{\textbf{\S \thebnomer.{ #1}}}\end{center}}

\newcommand{\thenv}[2]{%
\refstepcounter{snomer}
\par\addvspace{\medskipamount}\textbf{{#1} \thesnomer.}
{#2}\par\addvspace{\medskipamount}}


\renewcommand{\refname}{References}

\date{}
\title{The field of invariants for the adjoint action of the
Borel group in the nilradical of a~parabolic subalgebra}
\author{Victoria Sevostyanova\thanks{Partially
supported by Israel Scientific Foundation grant 797/14.}}

\maketitle

\begin{center}
\parbox[b]{330pt}{\small\textsc{Abstract.} In this paper the
field of invariants for the adjoint action of the Borel group in the
nilradical of a parabolic subalgebra is studied. We construct the
set of $B$-invariant rational functions generating the field of
invariants.}
\end{center}

\sect{Introduction}

Let $G$ be the general linear group $\mathrm{GL}(n,K)$ over an
algebraically closed field $K$ of characteristic zero. Let $B$ ($N$,
respectively) be its Borel (maximal unipotent, respectively)
subgroup, which consists of upper triangular matrices with nonzero
(unit, respectively) elements on the diagonal. We fix a parabolic
subgroup $P\supset B$. Let $\mathfrak{p}$, $\mathfrak{b}$ and
$\mathfrak{n}$ be the Lie subalgebras in $\mathfrak{gl}(n,K)$
correspon\-ding to $P$, $B$ and $N$, respectively. We represent
$\mathfrak{p}=\mathfrak{r}\oplus\mathfrak{m}$ as the direct sum of
the nilradical $\mathfrak{m}$ and a block diagonal subalgebra
$\mathfrak{r}$ with sizes of blocks $(r_1,r_2,\ldots,r_u)$. The
subalgebra $\mathfrak{m}$ is invariant relative to the adjoint
action of the group $P$:
$$\mbox{for any }g\in P\mbox{ we have
}x\in\mathfrak{m}\mapsto\mathrm{Ad}_gx=gxg^{-1}.$$ Therefore
$\mathfrak{m}$ is invariant relative to the adjoint action of the
subgroups $B$ and $N$. We extend this action to the representation
in the algebra $K[\mathfrak{m}]$ and in the field $K(\mathfrak{m})$:
$$\mbox{for any }g\in P\mbox{ we have
}f(x)\in K[\mathfrak{m}]\mapsto f(\mathrm{Ad}_{g^{-1}}x).$$

The complete description of the field of invariants
$K(\mathfrak{m})^N$ for any parabolic subalgebra is a result
of~\cite{S1}. In this paper a notion of an extended base is
introduced. The extended base is a subset of the set of positive
roots such that elements of the extended base correspond to a set of
algebraically independent $N$-invariants. These invariants generate
the field of invariants $K(\mathfrak{m})^N$.

The aim of this paper is to to study the field of invariants
$K(\mathfrak{m})^B$. It continues the series of
works~\cite{PS},\cite{S1},\cite{S2},\cite{S3}. In the paper we
introduce an analog of extended base for $B$-action. We determine a
subset $\Psi$ of the extended base. Every root of $\Psi$ corresponds
to a rational function. We show that these rational functions are
$B$-invariant (Lemmas~\ref{A_invariant} and~\ref{B_invariant}) and
algebraically independent (Proposition~\ref{A-B_independ}) and
generate the field of invariants $K(\mathfrak{m})^B$
(Corollary~\ref{field_B-invariant}). We also construct a
representative of any $B$-orbit in general position
(Theorem~\ref{Exist_representative_for_B}).

Further in the papers~\cite{S2},\cite{S3} the structure of the
algebra of invariants $K[\mathfrak{m}]^N$ is considered. If the
sizes of diagonal blocks are $(2,k,2)$, $k>2$, or $(1,2,2,1)$, then
the invariants constructed on the extended base do not generate the
algebra of invariants and the algebra of invariants is not free.
Besides, the additional invariants in both cases are constructed,
which together with the main list of invariants constructed on the
extended base generate the algebra of invariants
$K[\mathfrak{m}]^N$. Also, the relations between these invariants
are provided. In~\cite{S3} it was shown that the algebra of
invariants $K[\mathfrak{m}]^N$ is finitely generated. Since
$B=T\ltimes N$, where $T$ is the reductive group of nondegenerate
diagonal matrices, the algebra of invariants
$K[\mathfrak{m}]^B=K\big[K[\mathfrak{m}]^N\big]^T$ is finitely
generated too. The structures of $K[\mathfrak{m}]^N$ and
$K[\mathfrak{m}]^B$ seem to be very mysterious and are a
considerable challenge.


\sect{Main statements and definitions}

We begin with definitions. Let
$\mathfrak{b}=\mathfrak{n}\oplus\mathfrak{h}$ be a triangular
decomposition. Let $\Delta$ be the root system relative to
$\mathfrak{h}$ and let $\Delta^{\!+}$ be the set of positive roots.
Let $\{\varepsilon_i\}_{i=1}^{n}$ be the standard basis of
$\mathbb{C}^n$. Every positive root $\gamma$ in $\mathfrak{gl}(n,K)$
can be represented as $\gamma=\varepsilon_i-\varepsilon_j$,
$1\leqslant i<j\leqslant n$ (see \cite{GG}). We identify a root
$\gamma$ with the pair $(i,j)$ and the set of the positive roots
$\Delta^{\!+}$ with the set of pairs $(i,j)$, $i<j$. The system of
positive roots $\Delta^{\!+}_\mathfrak{r}$ of the reductive
subalgebra $\mathfrak{r}$ is a subsystem in $\Delta^{\!+}$.

Let $\{E_{i,j}:~i<j\}$ be the standard basis in $\mathfrak{n}$. Let
$E_\gamma$ denote the basis element $E_{i,j}$, where $\gamma=(i,j)$.

Let $M$ be a subset of $\Delta^{\!+}$ corresponding to
$\mathfrak{m}$ that is $$\mathfrak{m}=\bigoplus_{\gamma\in
M}E_{\gamma}.$$ We identify the algebra $K[\mathfrak{m}]$ with the
polynomial algebra in variables $x_{i,j}$, $(i,j)\in M$.

We define a relation in $\Delta^{\!+}$ setting $\gamma'>\gamma$
whenever $\gamma'-\gamma\in\Delta^{\!+}$. Note that the relation $>$
is not an order relation.

The roots $\gamma$ and $\gamma'$ are called \emph{comparable}, if
either $\gamma'>\gamma$ or $\gamma>\gamma'$.

We will introduce a subset $S$ in the set of positive roots such
that every root from this subset corresponds to some $N$-invariant.

\thenv{Definition}{A subset $S$ in $M$ is called a \emph{base} if
the elements in $S$ are not pairwise comparable and for any
$\gamma\in M\setminus S$ there exists $\xi\in S$ such that
$\gamma>\xi$.}

Let us show that the base exists. We need the following

\thenv{Definition\label{Def-minimal}}{Let $A$ be a subset in $M$. We
say that $\gamma$ is a \emph{minimal element} in $A$ if there is no
$\xi\in A$ such that $\xi<\gamma$.}

For a given parabolic subgroup we will construct a diagram in the
form of a square array. The cell of the diagram corresponding to a
root of $S$ is labeled by the symbol $\otimes$. Symbols $\times$
will be explained below.

\thenv{Example\label{Ex-2132}}{Diagram 1 represents the parabolic
subalgebra with sizes of its diagonal blocks $(2,1,3,2)$. In this
case minimal elements in $M$ are $(2,3)$, $(3,4)$ and $(6,7)$.}
\begin{center}\refstepcounter{diagram}
{\begin{tabular}{|p{0.1cm}|p{0.1cm}|p{0.1cm}|p{0.1cm}|p{0.1cm}|p{0.1cm}|p{0.1cm}|p{0.1cm}|c}
\multicolumn{2}{l}{{\small 1\quad2\!\!}}&\multicolumn{2}{l}{{\small
3\quad4\!\!}}&\multicolumn{2}{l}{{\small 5\quad 6\!\!}}&
\multicolumn{2}{l}{{\small 7\quad 8\!\!}}\\
\cline{1-8} \multicolumn{2}{|l|}{1}&&&$\!\otimes$&&&&{\small 1}\\
\cline{3-8} \multicolumn{2}{|r|}{1}&$\!\otimes$&&&&&&{\small 2}\\
\cline{1-8} \multicolumn{2}{|c|}{}&1&$\!\otimes$&&&&&{\small 3}\\
\cline{3-8} \multicolumn{3}{|c|}{}&\multicolumn{3}{|l|}{1}&$\!\times$&$\!\times$&{\small 4}\\
\cline{7-8} \multicolumn{3}{|c|}{}&\multicolumn{3}{|c|}{1}&$\!\times$&$\!\otimes$&{\small 5}\\
\cline{7-8} \multicolumn{3}{|c|}{}&\multicolumn{3}{|r|}{1}&$\!\otimes$&&{\small 6}\\
\cline{4-8} \multicolumn{6}{|c|}{}&\multicolumn{2}{|l|}{1}&{\small 7}\\
\multicolumn{6}{|c|}{}&\multicolumn{2}{|r|}{1}&{\small 8}\\
\cline{1-8} \multicolumn{8}{c}{Diagram \arabic{diagram}}\\
\end{tabular}}
\end{center}

We construct the base $S$ by the following algorithm.

\textsc{Step 1.} Put $M_0=M$ and $i=1$. Let $S_1$ be the set of
minimal elements in $M_0$.

\textsc{Step 2.} Put $M_i=M_{i-1}\setminus\big\{S_i\cup\{\gamma\in
M_{i-1}:\exists\xi\in S_i, \xi<\gamma\}\big\}$. Let $S_i$ be the set
of minimal elements of $M_{i-1}$. Increase $i$ by 1 and repeat Step
2 until $M_i$ is empty.

Denote $S=S_1\cup S_2\cup\ldots$ The base $S$ is unique.

We have $S_1=\{(2,3),(3,4),(6,7)\}$ and $S_2=\{(1,5),(5,8)\}$ in
Example~\ref{Ex-2132}.

\medskip
Let $(r_1,r_2,\ldots,r_s)$ be the sizes of the diagonal blocks in
$\mathfrak{r}$. Put $$R_k=\displaystyle\sum_{i=1}^{k}r_i.$$

Let us present $N$-invariant corresponding to a root of the base.
Consider the formal matrix $\mathbb{X}$ of variables
$$\left(\mathbb{X}\right)_{i,j}=\left\{
\begin{array}{ll}
x_{i,j}&\mbox{if }(i,j)\in M;\\
0&\mbox{otherwise.}
\end{array}\right.$$
The matrix $\mathbb{X}$ can be represented as a block matrix
$$\mathbb{X}=\left(\begin{array}{ccccc}
0&X_{1,2}&X_{1,3}&\ldots&X_{1,s}\\
0&0&X_{2,3}&\ldots&X_{1,s}\\
\ldots&\ldots&\ldots&\ldots&\ldots\\
0&0&0&\ldots&X_{s-1,s}\\
0&0&0&\ldots&0\\
\end{array}\right),$$
where the size of $X_{i,j}$ is $r_i\times r_j$,
\begin{equation}
X_{i,j}=\left(\begin{array}{cccc}
x_{R_{i-1}+1,R_{j-1}+1}&x_{R_{i-1}+1,R_{j-1}+2}&\ldots&x_{R_{i-1}+1,R_{j}}\\
x_{R_{i-1}+1,R_{j-1}+2}&x_{R_{i-1}+2,R_{j-1}+2}&\ldots&x_{R_{i-1}+2,R_{j}}\\
\ldots&\ldots&\ldots&\ldots\\
x_{R_{i},R_{j-1}+1}&x_{R_{i},R_{j-1}+2}&\ldots&x_{R_{i},R_{j}}\\
\end{array}\right).\label{X_ij}
\end{equation}

\thenv{Lemma\label{Lemma1}}{\emph{The roots corresponding to the
antidiagonal elements in} $X_{i,i+1}$ (\emph{from the lower left
element towards right upper direction}) \emph{are in the base}.}

Thus the roots of the base in the blocks $X_{i,i+1}$ are as follows.
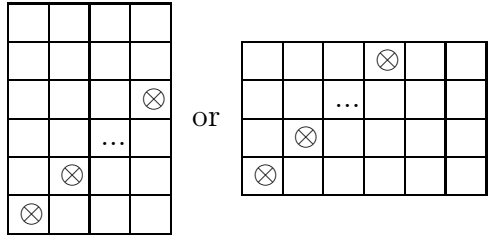
\begin{center}
{\begin{tabular}{|p{0.1cm}|p{0.1cm}|p{0.1cm}|p{0.1cm}|}
\cline{1-4} &&&\\
\cline{1-4} &&&\\
\cline{1-4} &&&$\!\otimes$\\
\cline{1-4} &&$\!...$&\\
\cline{1-4} &$\!\otimes$&&\\
\cline{1-4} $\!\otimes$&&&\\
\cline{1-4}
\end{tabular}\ \ \mbox{or }
\begin{tabular}{|p{0.1cm}|p{0.1cm}|p{0.1cm}|p{0.1cm}|p{0.1cm}|p{0.1cm}|}
\cline{1-6} &&&$\!\otimes$&&\\
\cline{1-6} &&$\!...$&&&\\
\cline{1-6} &$\!\otimes$&&&&\\
\cline{1-6} $\!\otimes$&&&&&\\
\cline{1-6}
\end{tabular}}
\end{center}

\textsc{Proof.} By definition~\ref{Def-minimal} for any $i$ the root
$(R_i,R_i+1)$ is minimal. Therefore $M\setminus M_1$ contains roots
corresponding to all cells in the row $R_i$ and the column $R_i+1$.
Hence $(R_i-1,R_i+2)\in S_2$ if $r_i,r_{i+1}>1$ and all roots of $M$
in the rows $R_i$, $R_i-1$ and in the columns $R_i+1$, $R_i+2$
belong to $M\setminus M_2$. Hence $(R_i-2,R_i+3)\in S_3$ if
$r_i,r_{i+1}>2$ etc.~$\Box$

There are roots in $S$ such that these roots do not correspond to
elements of the secondary diagonal in $X_{i,i+1}$, for example
$(1,5)$ in Example~\ref{Ex-2132}.

\thenv{Lemma\label{Lemma2}}{\emph{Let the maximal size of blocks in
$\mathfrak{r}$ between the $i$th and $j$th blocks be $k$ and $k<r_i$
and} $k<r_j$, \emph{then the roots of the base in the block
$X_{i,j}$ are following}
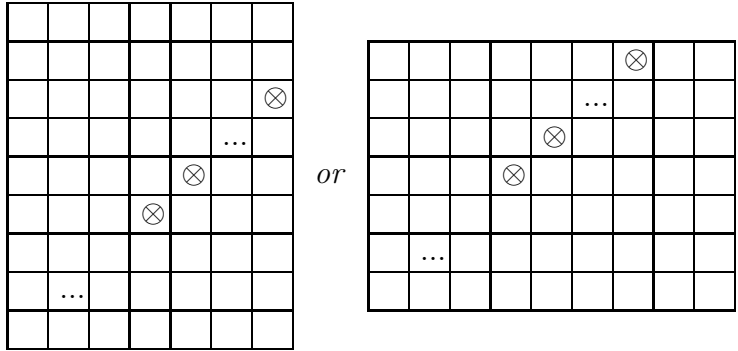
\begin{center}
{\begin{tabular}{|p{0.1cm}|p{0.1cm}|p{0.1cm}|p{0.1cm}|p{0.1cm}|p{0.1cm}|p{0.1cm}|}
\cline{1-7} &&&&&&\\
\cline{1-7} &&&&&&\\
\cline{1-7} &&&&&&$\!\otimes$\\
\cline{1-7} &&&&&$\!...$&\\
\cline{1-7} &&&&$\!\otimes$&&\\
\cline{1-7} &&&$\!\otimes$&&&\\
\cline{1-7} &&&&&&\\
\cline{1-7} &$\!...$&&&&&\\
\cline{1-7} &&&&&&\\
\cline{1-7}
\end{tabular}\ \mbox{ \emph{or} }
\begin{tabular}{|p{0.1cm}|p{0.1cm}|p{0.1cm}|p{0.1cm}|p{0.1cm}|
p{0.1cm}|p{0.1cm}|p{0.1cm}|p{0.1cm}|}
\cline{1-9} &&&&&&$\!\otimes$&&\\
\cline{1-9} &&&&&$\!...$&&&\\
\cline{1-9} &&&&$\!\otimes$&&&&\\
\cline{1-9} &&&$\!\otimes$&&&&&\\
\cline{1-9} &&&&&&&&\\
\cline{1-9} &$\!...$&&&&&&&\\
\cline{1-9} &&&&&&&&\\
\cline{1-9}
\end{tabular}}
\end{center}
\emph{where there are $k$ empty cells on the antidiagonal}
(\emph{from the lower left element towards right upper direction}).
\emph{There is a root of $S$ under every cell of these $k$ empty
cells and there is a root of $S$ to the left of every cell of these
$k$ empty cells}. \emph{Besides}, \emph{for the first diagram there
is not a root of $S$ to the left of every first empty rows and for
the second diagram there is not a root of $S$ under every last empty
columns}.}

\thenv{Example}{In Example~\ref{Ex-2132} the block $X_{1,3}$
corresponds to the following diagram.
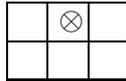
\begin{center}
{\begin{tabular}{|p{0.1cm}|p{0.1cm}|p{0.1cm}|}
\cline{1-3} &$\!\otimes$&\\
\cline{1-3} &&\\
\cline{1-3}
\end{tabular}}
\end{center}}

\textsc{Proof.} We show the statement by induction on $k$. To prove
the base of induction we consider $k=0$ and $j=i+1$, the case
follows from Lemma~\ref{Lemma1}. Suppose that the statement is true
for any maximal size of blocks between blocks $i$ and $j$ less that
$k$. Let us prove for $k$. Let $k$ is the size of $m$th block,
$k=r_m$. Consider the block $X_{i,m}$. The sizes of blocks between
$i$th and $m$th blocks are less than $k$. By assumption of
induction, the structure of $X_{i,m}$ is as in the first diagram in
Lemma~\ref{Lemma2}. Therefore there is a root of $S$ to the left of
every last $k$ rows and there is not a root of $S$ to the left of
the other rows of $X_{i,m}$. Hence the same is true for $X_{i,j}$.

Similarly, there is a root of $S$ under every first $k$ columns and
there is not a root of $S$ under the other columns of $X_{i,j}$.
Hence we have the statement of the lemma.~$\Box$

\medskip
For any root $\gamma=(a,b)\in M$ let $S_\gamma=\{(i,j)\in
S:i>a,j<b\}$. Let $S_\gamma=\{(i_1,j_1),\ldots,(i_k,j_k)\}$. Note
that if $\gamma$ is minimal in $M$, then $S_{\gamma}=\emptyset$.
Denote by $M_\gamma$ a minor $\mathbb{X}_I^J$ of the matrix
$\mathbb{X}$ with ordered systems of rows $I$ and columns $J$, where
$$I=\mathrm{ord}\{a,i_1,\ldots,i_k\},\quad J=\mathrm{ord}\{j_1,\ldots,j_k, b\}.$$

\thenv{Example}{Let us continue Example~\ref{Ex-2132}. For the root
$(1,6)$ we have $S_{(1,6)}=\{(2,3),(3,4)\}$, $I=\{1,2,3\}$,
$J=\{3,4,6\}$, and
$$M_{(1,6)}=\left|\begin{array}{ccc}
x_{1,3}&x_{1,4}&x_{1,6}\\
x_{2,3}&x_{2,4}&x_{2,6}\\
0&x_{3,4}&x_{3,6}\\
\end{array}\right|.$$
All minors $M_{\xi}$ for $\xi\in S$ are following
$$M_{(2,3)}=x_{2,3},\ M_{(3,4)}=x_{3,4},\ M_{(6,7)}=x_{6,7},$$$$
M_{(5,8)}=\left|\begin{array}{cc}
x_{5,7}&x_{5,8}\\
x_{6,7}&x_{6,8}\\
\end{array}\right|,\
M_{(1,5)}=\left|\begin{array}{ccc}
x_{1,3}&x_{1,4}&x_{1,5}\\
x_{2,3}&x_{2,4}&x_{2,5}\\
0&x_{3,4}&x_{3,5}\\
\end{array}\right|.$$}

\thenv{Lemma\label{L-M-inv}}{\emph{For any $\xi\in S$ the minor
$M_{\xi}$ is $N$-invariant}.}

\thenv{Notation\label{Notation}}{The group $N$ is generated by the
one-parameter subgroups
$$g_{i,j}(t)=I+tE_{i,j},\mbox{ where }1\leqslant i<j\leqslant n$$
and $I$ is the identity matrix. The adjoint action of any
$g_{i,j}(t)$ makes the following transformations of a matrix:
\begin{itemize}
\item[1)] the $j$th row multiplied by $t$ is added to the $i$th row,
\item[2)] the $i$th column multiplied by $-t$ is added to the $j$th
column, i.e. for a variable $x_{a,b}$ we have
$$\mathrm{Ad}_{g_{i,j}^{-1}(t)}x_{a,b}=\left\{\begin{array}{ll}
x_{a,b}+tx_{j,b}&\mbox{if }a=i;\\
x_{a,b}-tx_{a,i}&\mbox{if }b=j;\\
x_{a,b}&\mbox{otherwise}.
\end{array}\right.$$
\end{itemize}}

\textsc{Proof.} By the notation it is sufficient to prove that for
any $\xi=(k,m)\in S$ the minor $M_{\xi}$ is invariant under the
adjoint action of $g_{i,j}(t)$ for any $i<j$. If $i<k$, then the
$i$th row does not belong to the minor $M_{\xi}$ and the adding of
the $j$th row to the $i$th row leaves $M_{\xi}$ unchanged. Let
$M_{\xi}=\mathbb{X}_I^J$ for some collections of rows $I$ and
columns $J$. If $i\geqslant k$, then since the numbers in $I$ are
consecutive, the number of any nonzero row $j$ at the intersection
with columns $J$ belongs to $I$. Then the adding of the $j$th row to
the $i$th row leaves $M_{\xi}$ unchanged again. Using the similar
reasoning for columns, we get that $M_{\xi}$ is
$N$-invariant.~$\Box$

\medskip
The set $\{M_{\xi},\ \xi\in S\}$ does not generate all the
$N$-invariants. There is the other series of $N$-invariants. To
present it we need

\thenv{Definition}{An ordered set of positive roots
$$\{\varepsilon_{i_1}-\varepsilon_{j_1},
\varepsilon_{i_2}-\varepsilon_{j_2},\ldots,
\varepsilon_{i_s}-\varepsilon_{j_s}\}$$ is called a \emph{chain} if
$j_1=i_2,j_2=i_3,\ldots,j_{s-1}=i_s$.}

\thenv{Definition}{We say that two roots $\xi,\xi'\in S$ form an
\emph{admissible pair} $q=(\xi,\xi')$ if there exists $\alpha_q$ in
the set $\Delta^{\!+}_\mathfrak{r}$ corresponding to the reductive
part $\mathfrak{r}$ such that the ordered set of roots
$\{\xi,\alpha_q,\xi'\}$ is a chain. In other words, roots
$\xi=\varepsilon_i-\varepsilon_j$ and
$\xi'=\varepsilon_k-\varepsilon_l$ are an admissible pair if
$\alpha_q=\varepsilon_j-\varepsilon_k\in\Delta^{\!+}_\mathfrak{r}$.
Note that the root $\alpha_q$ is uniquely determined by $q$.}

\thenv{Example\label{Ex-2132-admissible_pair}}{In the case of
Diagram 1 we have three admissible pairs $q_1=(\xi_1,\xi_3)$,
$q_2=(\xi_2,\xi_3)$, $q_3=(\xi_1,\xi_4)$, where $\xi_1=(2,3)$,
$\xi_2=(1,5)$, $\xi_3=(6,7)$, and $\xi_4=(5,8)$.}

Let the set $Q:=Q(\mathfrak{p})$ consist of admissible pairs. For
every admissible pair $q=(\xi,\xi')$ we construct a positive root
$\varphi_q=\alpha_q+\xi'$, where $\{\xi,\alpha_q,\xi'\}$ is a chain.
Consider the subset $\Phi=\{\varphi_q:~q\in Q\}$ in the set of
positive roots. The cell of the diagram corresponding to a root of
$\Phi$ is labeled by $\times$.

\thenv{Example}{The roots of $\Phi$ for the admissible pairs in
Example~\ref{Ex-2132-admissible_pair} are $\varphi_{q_1}=(4,7)$,
$\varphi_{q_2}=(5,7)$, $\varphi_{q_3}=(4,8)$.}

Now we are ready to present the $N$-invariant corresponding to a
root $\varphi\in\Phi$.

Let admissible pair $q=(\xi,\xi')$ correspond to $\varphi_q\in\Phi$.
We construct the polynomial
\begin{equation}
L_{\varphi_q}=\sum_{\scriptstyle\alpha_1,\alpha_2\in\Delta^{\!+}_\mathfrak{r}\cup\{0\}
\atop\scriptstyle\alpha_1+\alpha_2=\alpha_q}
M_{\xi+\alpha_1}M_{\alpha_2+\xi'}.\label{L_q}
\end{equation}

\thenv{Example}{Continuing the previous example, we have
$$L_{(4,7)}=x_{3,4}x_{4,7}+x_{3,5}x_{5,7}+x_{3,6}x_{6,7},$$
$$L_{(4,8)}=x_{3,4}\left|\begin{array}{cc}
x_{4,7}&x_{4,8}\\
x_{6,7}&x_{6,8}\\
\end{array}\right|+x_{3,5}\left|\begin{array}{cc}
x_{5,7}&x_{5,8}\\
x_{6,7}&x_{6,8}\\
\end{array}\right|,$$
$$L_{(5,7)}=\left|\begin{array}{ccc}
x_{1,3}&x_{1,4}&x_{1,5}\\
x_{2,3}&x_{2,4}&x_{2,5}\\
0&x_{3,4}&x_{3,5}\\
\end{array}\right|x_{5,7}+\left|\begin{array}{ccc}
x_{1,3}&x_{1,4}&x_{1,6}\\
x_{2,3}&x_{2,4}&x_{2,6}\\
0&x_{3,4}&x_{3,6}\\
\end{array}\right|x_{6,7}.$$}

\thenv{Lemma}{\emph{The polynomial $L_{\varphi}$ is $N$-invariant
for any} $\varphi=\varphi_q\in\Phi$, $q=(\xi,\xi')$.}

\textsc{Proof.} By Notation~\ref{Notation} it is sufficient to prove
for the action of $g_{i,j}(t)$. Let $\xi=(a,b)$,~$\xi'=(a',b')$.
Using the definition of admissible pair, we have $a<b<a'<b'$,
$\alpha_q=(b,a')\in\Delta_{\mathfrak{r}}^{\!+}$, and
$\varphi=(b,b')$. If $i<b$ or $j>a'$, then using the same arguments
as in the proof of the invariance of $M_{\xi}$ for $\xi\in S$, we
have that the minors of the right part of (\ref{L_q}) are
$g_{i,j}(t)$-invariant.

Let $b\leqslant i<j\leqslant a'$. Denote $\gamma_1=(b,i)$,
$\gamma_2=(j,a')$, $\beta=(i,j)$; then
$\alpha_q=\gamma_1+\beta+\gamma_2$ and
$\gamma_1+\beta,\beta+\gamma_2\in\Delta_{\mathfrak{r}}^{\!+}\cup\{0\}$.
We have
\begin{equation}
\left\{\begin{array}{l} T_{g_{i,j}(t)}M_{\xi+\gamma_1+\beta}=
M_{\xi+\gamma_1+\beta}+tM_{\xi+\gamma_1},\\
T_{g_{i,j}(t)}M_{\beta+\gamma_2+\xi'}=
T_{\beta+\gamma_2+\xi'}-tM_{\gamma_2+\xi'}.\\
\end{array}\right.\label{T_g(M_xi)}
\end{equation}
The other minors of (\ref{L_q}) are invariant under the action of
$g_{i,j}(t)$. Combining (\ref{L_q}) and (\ref{T_g(M_xi)}), we get
$$\left(T_{g_{i,j}(t)}L_{\varphi}\right)-L_{\varphi}=M_{\xi+\gamma_1}
\left(M_{\beta+\gamma_2+\xi'}-t M_{\gamma_2+\xi'}\right)+$$
$$\left(M_{\xi+\gamma_1+\beta}+t
M_{\xi+\gamma_1}\right)M_{\gamma_2+\xi'}-
M_{\xi+\gamma_1}M_{\beta+\gamma_2+\xi'}-
M_{\xi+\gamma_1+\beta}M_{\gamma_2+\xi'}=0.~\Box$$

Thus we proved the first part of

\thenv{Theorem\label{M-L_independ}}{\emph{For an arbitrary parabolic
subalgebra}, \emph{the system of polynomials}
\begin{equation}
\{M_\xi,~\xi\in S,~L_{\varphi},~\varphi\in\Phi,\}\label{system-M-L}
\end{equation}
\emph{is contained in $K[\mathfrak{m}]^N$ and is algebraically
independent over $K$}.}

To show the algebraic independence, consider the restriction
homo\-mor\-phism $f\mapsto f|_\mathcal{Y}$, where
$$\mathcal{Y}=\left\{\sum_{\xi\in S\cup\Phi}c_{\xi}E_\xi:
c_{\xi}\neq0\ \forall\xi\in S\cup\Phi\right\},$$ from
$K[\mathfrak{m}]$ to the polynomial algebra $K[\mathcal{Y}]$ of
$x_\xi$, $\xi\in S$, and of $x_\varphi$, $\varphi\in\Phi$. Direct
calculations show that the system of the images
$$\left\{M_\xi|_{\mathcal{Y}},\xi\in S,L_\varphi|_{\mathcal{Y}},
\varphi\in\Phi\right\}$$ is algebraically independent over $K$.
Therefore, the system~(\ref{system-M-L}) is algebraically
independent over $K$ (see details in~\cite{PS}).

\thenv{Definition}{The set $S\cup\Phi$ is called an \emph{extended
base}.}

\thenv{Definition}{The matrices of $\mathcal{Y}$ are called
\emph{canonical}.}

By~\cite{S1} one has the following theorems.

\thenv{Theorem\label{Exist_of_representative}}{\emph{There exists a
nonempty Zariski-open subset $W\subset\mathfrak{m}$ such that the
$N$-orbit of any $x\in W$ intersects $\mathcal{Y}$ at a unique
point}.}

\thenv{Theorem\label{Th_invariant_field}}{\emph{The field of
invariants $K(\mathfrak{m})^N$ is the field of rational functions
of} $M_\xi$, $\xi\in S$, \emph{and} $L_{\varphi}$,
$\varphi\in\Phi$.}


\sect{Invariants of the Borel group in the nilradical of a~parabolic
subalgebra}

In this paragraph we will introduce some subset $\Psi\subset
S\cup\Phi$ and construct a correspon\-ding invariant in
$K[\mathfrak{m}]^B$ to every root $\xi\in \Psi$.

We describe two series of invariants in $K[\mathfrak{m}]^B$. The
first (second, respectively) series corresponds to a subset
$\Psi_1\subset\Psi$ ($\Psi_2\subset\Psi$, respectively) and
$\Psi=\Psi_1\sqcup\Psi_2$. The first series is following. Let
$\psi=(i,b)$ be a root in $\Phi$ such that there exist three roots
$\xi_1=(i,a)$, $\xi_2=(j,a)$, $\xi_3=(j,b)$ in $S\cup\Phi$, where
$i<j$ and $a<b$. Since there is no more than one root from $S$ in a
row and in a column, the roots $\xi_1,\xi_2$ belong to the set
$\Phi$. In other words, these four roots correspond to one of two
following diagrams
\begin{center}
{\begin{tabular} {p{0.1cm}p{0.1cm}p{0.1cm}|l}
\small{a}&\small{b}\\\cline{1-3}
\multicolumn{3}{c|}{}\\
$\times$&$\times$&&\small{i}\\
$\times$&$\otimes$&&\small{j}\\
\end{tabular}\qquad\mbox{or}\qquad
\begin{tabular}
{p{0.1cm}p{0.1cm}p{0.1cm}|l} \small{a}&\small{b}\\\cline{1-3}
\multicolumn{3}{c|}{}\\
$\times$&$\times$&&\small{i}\\
$\times$&$\times$&&\small{j}\\
\end{tabular}}
\end{center}
If there are such three roots $\xi_1,\xi_2,\xi_3$ with the above
described properties, then we say that $\psi$ belongs to the set
$\Psi_1$.

Since $\xi_2\in\Phi$, it corresponds to an admissible pair
$(\gamma,\gamma')$. Denote
$$A_{\psi}=\left\{
\begin{array}{ll}
\displaystyle\frac{L_{\psi}L_{\xi_2}}{L_{\xi_1}L_{\xi_3}}&\mbox{if
}\xi_3\in\Phi,\\
\\
\displaystyle\frac{L_{\psi}L_{\xi_2}}{L_{\xi_1}M_{\gamma}M_{\xi_3}}&\mbox{if
}\xi_3\in S.\\
\end{array}\right.$$ Evidently, the rational function $A_{\psi}$ is
$N$-invariant and depends on a choice of $\xi_1$, $\xi_2$, and
$\xi_3$. We will show below that $A_{\psi}$ is $B$-invariant.

\thenv{Example}{\label{Ex(2,3,2)}Consider the case when the diagonal
blocks of $\mathfrak{r}$ are $(2,3,2)$.
\begin{center}\refstepcounter{diagram}
{\begin{tabular}{|p{0.1cm}|p{0.1cm}|p{0.1cm}|p{0.1cm}|p{0.1cm}|p{0.1cm}|p{0.1cm}|c}
\multicolumn{2}{l}{{\small 1\quad 2}}&\multicolumn{2}{l}{{\small
3\quad 4}}&\multicolumn{2}{l}{{\small 5\quad 6}}&
\multicolumn{2}{l}{{\small 7}}\\
\cline{1-7} \multicolumn{2}{|l|}{1}&&$\!\otimes$&&&&{\small 1}\\
\cline{3-7} \multicolumn{2}{|r|}{1}&$\!\otimes$&&&&&{\small 2}\\
\cline{1-7} \multicolumn{2}{|c|}{}&\multicolumn{3}{|l|}{1}&$\!\times$&$\!\boxtimes$&{\small 3}\\
\cline{6-7} \multicolumn{2}{|c|}{}&\multicolumn{3}{|c|}{1}&$\!\times$&$\!\otimes$&{\small 4}\\
\cline{6-7} \multicolumn{2}{|c|}{}&\multicolumn{3}{|r|}{1}&$\!\otimes$&&{\small 5}\\
\cline{3-7} \multicolumn{5}{|c|}{}&\multicolumn{2}{|l|}{1}&{\small 6}\\
\multicolumn{5}{|c|}{}&\multicolumn{2}{|r|}{1}&{\small 7}\\
\cline{1-7} \multicolumn{7}{c}{Diagram \arabic{diagram}}\\
\end{tabular}}
\end{center}
The unique root $\psi\in\Psi_1$ is $(3,7)$. We will denote roots of
$\Psi$ by the symbol $\boxtimes$. In this case we have
$\xi_1=(3,6)$, $\xi_2=(4,6)$, and $\xi_3=(4,7)\in S$. The root
$\xi_3$ corresponds to the admissible pair $\big((1,4),(5,6)\big)$,
then
$$A_{(3,7)}=\frac{L_{(3,7)}L_{(4,6)}}{L_{(3,6)}M_{(1,4)}M_{(4,7)}},$$
where
$$L_{(3,7)}=x_{2,3}\left|\begin{array}{cc}
x_{3,6}&x_{3,7}\\
x_{5,6}&x_{5,7}\\
\end{array}\right|+x_{2,4}\left|\begin{array}{cc}
x_{4,6}&x_{4,7}\\
x_{5,6}&x_{5,7}\\
\end{array}\right|,$$
$$L_{(4,6)}=\left|\begin{array}{cc}
x_{1,3}&x_{1,4}\\
x_{2,3}&x_{2,4}\\
\end{array}\right|x_{4,6}+\left|\begin{array}{cc}
x_{1,3}&x_{1,5}\\
x_{2,3}&x_{2,5}\\
\end{array}\right|x_{5,6},$$
$$L_{(3,6)}=x_{2,3}x_{3,6}+
x_{2,4}x_{4,6}+x_{2,5}x_{5,6},$$ $$M_{(1,4)}=\left|\begin{array}{cc}
x_{1,3}&x_{1,4}\\
x_{2,3}&x_{2,4}\\
\end{array}\right|,\ M_{(4,7)}=\left|\begin{array}{cc}
x_{4,6}&x_{4,7}\\
x_{5,6}&x_{5,7}\\
\end{array}\right|.$$}

\medskip
If $\xi=(i,j)$, then denote $\mathrm{row}(\xi)=i$ and
$\mathrm{col}(\xi)=j$.

\thenv{Lemma\label{A_invariant}}{\emph{Rational functions}
$A_{\psi}$ , $\psi\in\Psi_1$, \emph{are invariant under the adjoint
action of} $B$.}

\textsc{Proof.} To prove the statement note that $A_{\psi}$ is
$N$-invariant. Since $B=T\ltimes N$, where $T$ is the set of
nondegenerate diagonal matrices, then it is sufficient to prove that
$A_{\psi}$ is $T$-invariant.

Let $\mathbb{X}_{I}^{J}$ be a square minor at the intersection of
rows $I=\{i_1,\ldots,i_k\}$ and columns $J=\{j_1,\ldots,j_k\}$. The
adjoint action of the diagonal element
$$t=\left(\begin{array}{ccc}
t_1&\ldots&0\\
\ldots&\ldots&\ldots\\
0&\ldots&t_n
\end{array}\right)=\mathrm{diag}(t_1,\ldots,t_n)$$
on $\mathbb{X}_{I}^{J}$ is follows
$$\mathrm{Ad}_t\,\mathbb{X}_{I}^{J}=\frac{t_I}{t_J}\mathbb{X}_{I}^{J},$$
where $t_I=t_{i_1}\ldots t_{i_k}$ and $t_J=t_{j_1}\ldots t_{j_k}$.

Take any $\varphi\in\Phi$. Let $(\mu_1,\mu_2)$ be an admissible pair
for $\varphi$. If $M_{\mu_1}=\mathbb{X}_{I}^{J_1}$ and
$M_{\mu_2}=\mathbb{X}_{I_2}^{J}$ for some collections of rows $I$
and $I_1$ and columns $J$ and $J_1$, then
\begin{equation}
D_{I,I'}^{J,J'}=\left|\begin{array}{cc}
\mathbb{X}_I^{J'}&(\mathbb{X}^2)_I^J\\
0&\mathbb{X}_{I'}^J\\
\end{array}\right|\label{L_phi-other}
\end{equation}
is called a combined minor, where
$I'=I_2\setminus\{\mathrm{row}(\mu_2)\}$ and
$J'=J_1\setminus\{\mathrm{col}(\mu_1)\}$. By Proposition~2.5
from~\cite{S2} the combined minor~(\ref{L_phi-other}) is equal to
the $N$-invariant $L_{\varphi}$.

\thenv{Example}{Let us show combined minors for invariants
$L_{\varphi}$ in Example~\ref{Ex(2,3,2)}. We have
$$L_{(3,6)}=D_{\{2\},\emptyset}^{\{6\},\emptyset}=
(\mathbb{X}^2)_{\{2\}}^{\{6\}},$$
$$L_{(3,7)}=D_{\{2\},\{5\}}^{\{6,7\},\emptyset}=
\left|\begin{array}{c}
(\mathbb{X}^2)_{\{2\}}^{\{6,7\}}\\
\mathbb{X}_{\{5\}}^{\{6,7\}}\\
\end{array}\right|,$$
$$L_{(4,6)}=D_{\{1,2\},\emptyset}^{\{6\},\{3\}}=
\left|\mathbb{X}_{1,2}^{3}\quad\!(\mathbb{X}^2)_{\{1,2\}}^{\{6\}}\right|.$$
Note that
$$D_{\{1,2\},\{5\}}^{\{6,7\},\{3\}}=\left|\begin{array}{cc}
\mathbb{X}_{\{1,2\}}^{\{3\}}&(\mathbb{X}^2)_{\{1,2\}}^{\{6,7\}}\\
0&\mathbb{X}_{\{5\}}^{\{6,7\}}\\
\end{array}\right|=M_{(2,3)}M_{(4,7)}.$$}

Since $\psi\in\Psi_1$, there are roots $\xi_1,\xi_2,\xi_3\in
S\cup\Phi$ such that $\psi=(i,b)$, $\xi_1=(i,a)$, $\xi_2=(j,a)$,
$\xi_3=(j,b)$ for some $i<j$ and $a<b$. Using the
submission~(\ref{L_phi-other}) we have
$$L_{\psi}=D_{I_1,I'}^{J_1,J'},\quad
L_{\xi_1}=D_{I_1,I''}^{J_2,J'},\quad
L_{\xi_2}=D_{I_2,I''}^{J_2,J''}$$ for some ordered sets of rows
$I_1\subset I_2$, $I'\supset I''$ and columns $J_1\supset J_2$,
$I_1\subset I_2$. If $\xi_3\in\Phi$, then
$L_{\xi_3}=D_{I_2,I'}^{J_1,J''}$. Let $\xi_3\in S$; since there
exists an admissible pair $(\gamma,\gamma')$ for $\xi_2$, then there
is the root $\gamma\in S$ in the column $j$. We have
$M_{\xi_3}M_{\gamma}=D_{I_2,I'}^{J_1,J''}$. In any case we get
$$A_{\psi}=\frac{D_{I_1,I'}^{J_1,J'}\cdot D_{I_2,I''}^{J_2,J''}}%
{D_{I_1,I''}^{J_2,J'}\cdot D_{I_2,I'}^{J_1,J''}}$$

Then
$$\mathrm{Ad}_t\, A_{\psi}=\frac{\displaystyle
\frac{t_{I_1}t_{I'}}{t_{J_1}t_{J'}}%
D_{I_1,I'}^{J_1,J'}\cdot\frac{t_{I_2}t_{I''}}{t_{J_2}t_{J''}}
D_{I_2,I''}^{J_2,J''}}{\displaystyle
\frac{t_{I_1}t_{I''}}{t_{J_2}t_{J'}}
D_{I_1,I''}^{J_2,J'}\cdot\frac{t_{I_2}t_{I'}}{t_{J_1}t_{J''}}
D_{I_2,I'}^{J_1,J''}}=A_{\psi},$$ i.e. a function $A_{\psi}$ is
$T$-invariant.~$\Box$

\medskip
Now we construct the second series of $B$-invariants.

Suppose for some $t$ we have $\psi=(R_{t-1}+k,R_t+1)\in\Phi$, where
$k<r_t$, i.e. $x_\psi$ is in the first column and in the $k$th row
of the block $X_{t-1,t}$. We say that the root $\psi$ belongs to the
set $\Psi_2$ of the second series if there exist $k$ roots
\begin{equation}
(R_{s-1}+1,j_1),(R_{s-1}+1,j_2),\ldots,(R_{s-1}+1,j_k)\in\Phi\label{k_roots}
\end{equation}
for some $j_1<j_2<\ldots<j_k$ and $s<t$. The roots~(\ref{k_roots})
correspond to variables in the first row of the blocks
$X_{s-1,s},X_{s-1,s+1},\ldots$ Suppose that the number $s$ is
maximal. We can assume without loss of generality that there are no
roots of $\Phi$ in the row $R_{s-1}+1$ between the
roots~(\ref{k_roots}), i.e. for any $a$ there is not a column $j$,
$j_a<j<j_{a+1}$, such that $(R_{s-1}+1,j)\in\Phi$. Denote by $\xi_1$
the last root $(R_{s-1}+1,j_k)$ in~(\ref{k_roots}). Let
$\gamma_1=(R_{s-1},R_{s-1}+1)$, evidently $\gamma_1\in S$. Since
$\psi\in\Phi$, then there is a corresponding admissible pair
$(\gamma_2,\gamma_3)$ for $\psi$, where $\gamma_2$ is in the
$(R_{t-1}+k)$th column and $\gamma_3=(R_t,R_t+1)$. See the diagram
below.

\begin{center}\refstepcounter{diagram}
{\begin{tabular}{p{0.1cm}|p{0.1cm}|p{0.1cm}|p{0.1cm}|p{0.1cm}|p{0.1cm}|
p{0.1cm}|p{0.1cm}|p{0.1cm}|p{0.1cm}|p{0.1cm}|p{0.1cm}|p{0.1cm}ll}
&$\!\gamma_1$&&&&&&&&&&&&&{\small $R_{s-1}$}\\
\cline{1-13} &\multicolumn{5}{|l|}{1}&&&$\!\xi_1$&&&$\!\gamma_5$&&&{\small $R_{s-1}+1$}\\
\cline{7-13} &\multicolumn{5}{|l|}{\quad\,...}&&&&&$\!...$&&\\
\cline{7-13} &\multicolumn{5}{|c|}{1}&&&&$\!\gamma_4$&&&\\
\cline{7-13} &\multicolumn{5}{|c|}{\qquad\,1}&&&$\!\gamma_2$&&&&\\
\cline{7-13} &\multicolumn{5}{|r|}{...\!}&&&&&&&\\
\cline{2-13} \multicolumn{6}{c|}{}&1&&&&&&&&{\small $R_{t-1}$}\\
\cline{7-13} \multicolumn{7}{c|}{}&\multicolumn{5}{|l|}{\!...}&\\
\cline{13-13} \multicolumn{7}{c|}{}&\multicolumn{5}{|l|}{\quad\,1}&$\!\psi$&&{\small $R_{t-1}+k$}\\
\cline{13-13} \multicolumn{7}{c|}{}&\multicolumn{5}{|c|}{1}&\\
\cline{13-13} \multicolumn{7}{c|}{}&\multicolumn{5}{|c|}{\qquad\,...}&\\
\cline{13-13} \multicolumn{7}{c|}{}&\multicolumn{5}{|r|}{1\!}&$\!\gamma_3$&&{\small $R_t$}\\
\cline{8-13} \multicolumn{12}{c|}{}&1\\
\multicolumn{13}{c}{Diagram \arabic{diagram}\label{diag1}}\\
\end{tabular}}
\end{center}

Let us show that the sizes of blocks between the $s$th and $t$th
blocks in $\mathfrak{r}$ are less or equal to $k$. Suppose $r_p>k$,
$s<p<t$. Therefore by Lemmas~\ref{Lemma1} and \ref{Lemma2} there are
$r_p$ roots of $S$ in rows $R_{p-1}+1,\ldots,R_p$, i.e. there is a
root in every row to the right of $p$th block. Every of these roots
of $S$ in rows $R_{p-1},\ldots,R_p$ and a root $(R_{p-1},R_{p-1}+1)$
make an admissible pair. Hence there exist $r_p-1$ roots of $\Phi$
in the row $R_{p-1}+1$. Since $r_p-1\geqslant k$, it is a
contradiction with the maximality of $s$.

From Lemmas \ref{Lemma1} and \ref{Lemma2} it follows that there
exists a root $\gamma_4\in S$ in the column $R_{t-1}+k+1$. Since the
sizes of blocks between the $s$th and $t$th blocks in $\mathfrak{r}$
are less or equal to $k$, $\gamma_4$ is to the right of the $s$th
block in $\mathfrak{r}$.

There exist three cases $r_s=r_t$, $r_s<r_t$, and $r_s>r_t$. We
illustrate all three cases and point all roots $\psi$, $\xi_1$,
$\gamma_i$ on Diagrams~\ref{diag1}, \ref{diag2}, and~\ref{diag3}.
Diagram 3 (4 and 5, respectively) corresponds to the case $r_s=r_t$
($r_s<r_t$ and $r_s>r_t$, respectively).

\begin{center}\refstepcounter{diagram}
{\begin{tabular}{p{0.1cm}|p{0.1cm}|p{0.1cm}|p{0.1cm}|p{0.1cm}|p{0.1cm}|
p{0.1cm}|p{0.1cm}|p{0.1cm}|p{0.1cm}|p{0.1cm}|p{0.1cm}|p{0.1cm}|p{0.1cm}ll}
&$\!\gamma_1$&&&&&&&&&&&&&&{\small $R_{s-1}$}\\
\cline{1-14} &\multicolumn{5}{|l|}{1}&&&$\!\xi_1$&&&$\!\gamma_5$&&&&{\small $R_{s-1}+1$}\\
\cline{7-14} &\multicolumn{5}{|l|}{\quad\,...}&&&&&$\!...$&&\\
\cline{7-14} &\multicolumn{5}{|c|}{1}&&&&$\!\gamma_4$&&&\\
\cline{7-14} &\multicolumn{5}{|c|}{\qquad\,1}&&&$\!\gamma_2$&&&&\\
\cline{7-14} &\multicolumn{5}{|r|}{...\!}&&&&&&&\\
\cline{2-14} \multicolumn{6}{c|}{}&1&&&&&&&&&{\small $R_{t-1}$}\\
\cline{7-14} \multicolumn{7}{c|}{}&\multicolumn{6}{|l|}{\!...}&\\
\cline{14-14} \multicolumn{7}{c|}{}&\multicolumn{6}{|l|}{\quad\,1}&$\!\psi$&&{\small $R_{t-1}+k$}\\
\cline{14-14} \multicolumn{7}{c|}{}&\multicolumn{6}{|c|}{1\quad\,}&\\
\cline{14-14} \multicolumn{7}{c|}{}&\multicolumn{6}{|c|}{\quad\,\,...}&\\
\cline{14-14} \multicolumn{7}{c|}{}&\multicolumn{6}{|c|}{\qquad\qquad1}&$\!\xi_2$\\
\cline{14-14} \multicolumn{7}{c|}{}&\multicolumn{6}{|r|}{1\!}&$\!\gamma_3$&&{\small $R_t$}\\
\cline{8-14} \multicolumn{13}{c|}{}&1\\
\multicolumn{14}{c}{Diagram \arabic{diagram}\label{diag2}}\\
\end{tabular}}
\end{center}

\begin{center}\refstepcounter{diagram}
{\begin{tabular}{p{0.1cm}|p{0.1cm}|p{0.1cm}|p{0.1cm}|p{0.1cm}|p{0.1cm}|
p{0.1cm}|p{0.1cm}|p{0.1cm}|p{0.1cm}|p{0.1cm}|p{0.1cm}|p{0.1cm}|p{0.1cm}ll}
&$\!\gamma_1$&&&&&&&&&&&&&&{\small $R_{s-1}$}\\
\cline{1-14} &\multicolumn{6}{|l|}{1}&&&$\!\xi_1$&&&$\!\xi_3$&&&{\small $R_{s-1}+1$}\\
\cline{8-14} &\multicolumn{6}{|l|}{\quad\,1}&&&&&&$\!\gamma_5$&\\
\cline{8-14} &\multicolumn{6}{|c|}{...\quad\,}&&&&&&$\!...$&\\
\cline{8-14} &\multicolumn{6}{|c|}{\quad\,\,1}&&&&$\!\gamma_4$&&&\\
\cline{8-14} &\multicolumn{6}{|c|}{\qquad\qquad1}&&&$\!\gamma_2$&&&&\\
\cline{8-14} &\multicolumn{6}{|r|}{...\!}&&&&&&&\\
\cline{2-14} \multicolumn{7}{c|}{}&1&&&&&&&&{\small $R_{t-1}$}\\
\cline{8-14} \multicolumn{8}{c|}{}&\multicolumn{5}{|l|}{\!...}&\\
\cline{14-14} \multicolumn{8}{c|}{}&\multicolumn{5}{|l|}{\quad\,1}&$\!\psi$&&{\small $R_{t-1}+k$}\\
\cline{14-14} \multicolumn{8}{c|}{}&\multicolumn{5}{|c|}{1}&\\
\cline{14-14} \multicolumn{8}{c|}{}&\multicolumn{5}{|c|}{\qquad\,...}&\\
\cline{14-14} \multicolumn{8}{c|}{}&\multicolumn{5}{|r|}{1\!}&$\!\gamma_3$&&{\small $R_t$}\\
\cline{9-14} \multicolumn{13}{c|}{}&1\\
\multicolumn{14}{c}{Diagram \arabic{diagram}\label{diag3}}\\
\end{tabular}}
\end{center}

Let us analyze these three cases.
\begin{enumerate}
\item If $r_s=r_t$, then denote $\gamma_5=(R_{s-1}+1,R_t)$.
Obviously $\gamma_5\in S$.
\item If $r_s<r_t$, then there exists a root in $S$ in the row
$R_{s-1}+1$. Denote this root by $\gamma_5$. In this case the pair
$(\gamma_5,\gamma_3)$ is admissible, therefore the pair corresponds
to  some $\xi_2\in\Phi$. The root $\xi_2$ is in the column $R_t+1$.
\item In the case $r_s>r_t$ there exists a root in $S$ in the column
$R_t$, we denote it by $\gamma_5$. Then $(\gamma_1,\gamma_5)$ is
admissible and we denote a corresponding root
$(R_{s-1}+1,R_t)\in\Phi$ by $\xi_3$.
\end{enumerate}

Denote
$$D=\left\{\begin{array}{ll}
M_{\gamma_1}M_{\gamma_3}M_{\gamma_5}&\mbox{if
}r_s=r_t;\\
M_{\gamma_1}L_{\xi_2}&\mbox{if }r_s<r_t;\\
M_{\gamma_3}L_{\xi_3}&\mbox{if }r_s>r_t.\\
\end{array}\right.$$
Now we are ready to show the second series of $B$-invariants.
Firstly, consider the simplest case, when $t=s+1$ and $r_s=2$ or
$r_t=2$ see the diagrams below.
\begin{center}
{\refstepcounter{diagram}
\begin{tabular}{p{0.1cm}|p{0.1cm}|p{0.1cm}|p{0.1cm}|p{0.1cm}|p{0.1cm}}
1&$\!\gamma_1$&&&&\\
\cline{1-6} &\multicolumn{2}{|l|}{1}&$\!\xi_1$&$\!\gamma_4$&\\
\cline{4-6} &\multicolumn{2}{|r|}{1}&$\!\gamma_2$&&\\
\cline{2-6} \multicolumn{3}{c|}{}&\multicolumn{2}{|l|}{1}&$\!\psi$\\
\cline{6-6} \multicolumn{3}{c|}{}&\multicolumn{2}{|r|}{1}&$\!\gamma_3$\\
\cline{4-6} \multicolumn{5}{c|}{}&1\\
\multicolumn{6}{c}{Diagram \arabic{diagram}}\\
\end{tabular}\quad\refstepcounter{diagram}
\begin{tabular}{p{0.1cm}|p{0.1cm}|p{0.1cm}|p{0.1cm}|p{0.1cm}|p{0.1cm}|p{0.1cm}|p{0.1cm}}
1&$\!\gamma_1$&&&&&&\\
\cline{1-8} &\multicolumn{4}{|l|}{1}&$\!\xi_1$&$\!\xi_3$&\\
\cline{6-8} &\multicolumn{4}{|l|}{\quad...}&$\!...$&$\!...$&\\
\cline{6-8} &\multicolumn{4}{|c|}{\quad1}&&$\!\gamma_4$&\\
\cline{6-8} &\multicolumn{4}{|r|}{1}&$\!\gamma_2$&&\\
\cline{2-8} \multicolumn{5}{l|}{}&\multicolumn{2}{|l|}{1}&$\!\psi$\\
\cline{8-8} \multicolumn{5}{r|}{}&\multicolumn{2}{|r|}{1}&$\!\gamma_3$\\
\cline{6-8} \multicolumn{7}{c|}{}&1\\
\multicolumn{8}{c}{Diagram \arabic{diagram}}\\
\end{tabular}\quad\refstepcounter{diagram}
\begin{tabular}{p{0.1cm}|p{0.1cm}|p{0.1cm}|p{0.1cm}|p{0.1cm}|p{0.1cm}|p{0.1cm}|p{0.1cm}}
1&$\!\gamma_1$&&&&&&\\
\cline{1-8} &\multicolumn{2}{|l|}{1}&$\!\xi_1$&$\!\gamma_4$&&&\\
\cline{4-8} &\multicolumn{2}{|r|}{1}&$\!\gamma_2$&&&&\\
\cline{2-8} \multicolumn{3}{c|}{}&\multicolumn{4}{|l|}{1}&$\!\psi$\\
\cline{8-8} \multicolumn{3}{r|}{}&\multicolumn{4}{|l|}{\quad1}&$\!\xi_2$\\
\cline{8-8} \multicolumn{3}{l|}{}&\multicolumn{4}{|c|}{\quad...}&$\!...$\\
\cline{8-8} \multicolumn{3}{r|}{}&\multicolumn{4}{|r|}{1}&$\!\gamma_3$\\
\cline{4-8} \multicolumn{7}{c|}{}&1\\
\multicolumn{8}{c}{Diagram \arabic{diagram}}\\
\end{tabular}}
\end{center}
In these cases denote
$$B_{\psi}=\frac{L_{\psi}L_{\xi_1}}{D}.$$

\thenv{Example\label{Ex-1221}}{The simplest case $(1,2,2,1)$ is
following.
\begin{center}\refstepcounter{diagram}
{\begin{tabular}{|p{0.1cm}|p{0.1cm}|p{0.1cm}|p{0.1cm}|p{0.1cm}|p{0.1cm}|c}
\multicolumn{2}{l}{{\small 1\quad 2}}&\multicolumn{2}{l}{{\small
3\quad 4}}&\multicolumn{2}{l}{{\small 5\quad 6}}\\
\cline{1-6} 1&$\!\otimes$&&&&&{\small 1}\\
\cline{1-6} &\multicolumn{2}{|l|}{1}&$\!\times$&$\!\otimes$&&{\small 2}\\
\cline{4-6} &\multicolumn{2}{|r|}{1}&$\otimes$&&&{\small 3}\\
\cline{2-6} \multicolumn{3}{|c|}{}&\multicolumn{2}{|l|}{1}&$\!\boxtimes$&{\small 4}\\
\cline{6-6} \multicolumn{3}{|c|}{}&\multicolumn{2}{|r|}{1}&$\!\otimes$&{\small 5}\\
\cline{4-6} \multicolumn{5}{|c|}{}&1&{\small 6}\\
\cline{1-6} \multicolumn{7}{c}{Diagram \arabic{diagram}}\\
\end{tabular}}
\end{center}
We have a unique root $\psi=(4,6)\in\Psi_2$ in $\Psi$ and the unique
$B$-invariant.
$$B_{\psi}=\frac{L_{(2,4)}L_{(4,6)}}{M_{(1,2)}M_{(5,6)}M_{(2,5)}}=
\frac{(x_{1,2}x_{2,4}+x_{1,3}x_{3,4})(x_{3,4}x_{4,6}+x_{3,5}x_{5,6})}%
{\displaystyle x_{1,2}x_{5,6}\left|\begin{array}{cc}
x_{2,4}&x_{2,5}\\
x_{3,4}&x_{3,5}\\
\end{array}\right|}.$$}

Let us construct invariants for more difficult case. We define a
relation in $\Delta^{\!+}$ such that for roots $\gamma=(a_1,b_1)$
and $\gamma'=(a_2,b_2)$ we have $\gamma\succ\gamma'$ whenever
$a_1<a_2$ and $b_1>b_2$.

If $r_s>2$ and $r_t>2$ or $s\neq t+1$ and $r_s,r_t\geqslant2$, then
a $B$-invariant of the second series is defined as follows:
$$B_{\psi}=\frac{\displaystyle L_{\xi_1}L_{\psi}\cdot
\prod_{\mu\prec\gamma_5}M_{\mu}}{\displaystyle
D\cdot\prod_{\mu\prec\gamma_4}M_{\mu}\cdot
\prod_{\mu'\prec\mu\prec\gamma_4}M_{\mu'}}.$$ The product in the
numerator (the first product in the denominator resp.) are taken on
all roots $\mu\prec\gamma_5$ ($\mu\prec\gamma_4$ resp.) such that
$\mu\in S$ and $\mu$ is maximal in the sense of the relation
$\prec$. For example,
$\displaystyle\prod_{\mu\prec(1,5)}M_{\mu}=M_{(2,3)}M_{(3,4)}$ for
Diagram 1. The second product in the denominator is taken on all
roots $\mu'\prec\mu$, where $\mu\prec\gamma_4$, such that
$\mu,\mu'\in S$ and $\mu,\mu'$ are maximal in the sense of the
relation $\prec$.

\thenv{Example}{Consider more difficult example than
Example~\ref{Ex-1221}. Let $(2,2,3,3,2)$ be the sizes of diagonal
blocks in $\mathfrak{r}$.
\begin{center}\refstepcounter{diagram}
{\begin{tabular}{|p{0.1cm}|p{0.1cm}|p{0.1cm}|p{0.1cm}|p{0.1cm}|p{0.1cm}|
p{0.1cm}|p{0.1cm}|p{0.1cm}|p{0.1cm}|p{0.1cm}|p{0.1cm}|c}
\multicolumn{2}{l}{{\small 1\!\quad2}}&\multicolumn{2}{l}{{\small
3\!\quad4}}&\multicolumn{2}{l}{{\small 5\!\quad 6}}&
\multicolumn{2}{l}{{\small 7\!\quad 8}}&\multicolumn{2}{l}{{\small
9\ \ 10\!\!\!}}&\multicolumn{2}{l}{{\small \!11\ \,12\!\!\!\!\!}}\\
\cline{1-12} \multicolumn{2}{|l|}{1}&&$\!\otimes$&&&&&&&&&{\small 1}\\
\cline{3-12} \multicolumn{2}{|r|}{1}&$\!\otimes$&&&&&&&&&&{\small 2}\\
\cline{1-12} \multicolumn{2}{|c|}{}&\multicolumn{2}{|l|}{1}&$\!\times$&$\!\otimes$&&&&&&&{\small 3}\\
\cline{5-12} \multicolumn{2}{|c|}{}&\multicolumn{2}{|r|}{1}&$\!\otimes$&&&&&&&&{\small 4}\\
\cline{3-12} \multicolumn{4}{|c|}{}&\multicolumn{3}{|l|}{1}&$\!\boxtimes$&$\!\boxtimes$&$\!\otimes$&&&{\small 5}\\
\cline{8-12} \multicolumn{4}{|c|}{}&\multicolumn{3}{|c|}{1}&$\!\times$&$\!\otimes$&&&&{\small 6}\\
\cline{8-12} \multicolumn{4}{|c|}{}&\multicolumn{3}{|r|}{1}&$\!\otimes$&&&&&{\small 7}\\
\cline{5-12} \multicolumn{7}{|c|}{}&\multicolumn{3}{|l|}{1}&$\!\boxtimes$&$\!\boxtimes$&{\small 8}\\
\cline{11-12} \multicolumn{7}{|c|}{}&\multicolumn{3}{|c|}{1}&$\!\boxtimes$&$\!\otimes$&{\small 9}\\
\cline{11-12} \multicolumn{7}{|c|}{}&\multicolumn{3}{|r|}{1}&$\!\otimes$&&{\small 10}\\
\cline{8-12} \multicolumn{10}{|c|}{}&\multicolumn{2}{|l|}{1}&{\small 11}\\
\multicolumn{10}{|c|}{}&\multicolumn{2}{|r|}{1}&{\small 12}\\
\cline{1-12} \multicolumn{12}{c}{Diagram \arabic{diagram}}\\
\end{tabular}}
\end{center}
We have $\Psi_2=\big\{(5,8),(8,11),(9,11)\big\}$. For the root
$(5,8)$ we have the simple case of the second series
$$B_{(5,8)}=\frac{L_{(5,8)}L_{(3,5)}}{L_{(6,8)}M_{(2,3)}}.$$ For
roots $(8,11)$ and $(9,11)$ one has $\gamma_1=(4,5)$,
$\gamma_3=(10,11)$, $\gamma_5=(5,10)$, and
$\displaystyle\prod_{\mu\prec\gamma_5}M_{\mu}=M_{(6,9)}$. For the
root $(8,11)$ we have $\gamma_4=(6,9)$,
$\displaystyle\prod_{\mu\prec\gamma_4}M_{\mu}=M_{(7,8)}$, and there
is not a root $\mu'$ such that $\mu'\prec\mu\prec\gamma_4$, then
$$B_{(8,11)}=\frac{L_{(8,11)}L_{(5,8)}M_{(6,9)}}%
{M_{(4,5)}M_{(5,10)}M_{(10,11)}M_{(7,8)}}.$$ For $(9,11)$ we get
$\gamma_4=(5,10)$,
$\displaystyle\prod_{\mu\prec\gamma_4}M_{\mu}=M_{(6,9)}$,
$\displaystyle\prod_{\mu'\prec\mu\prec\gamma_4}M_{\mu'}=M_{(7,8)}$,
and
$$B_{(9,11)}=\frac{L_{(9,11)}L_{(5,9)}}%
{M_{(4,5)}M_{(5,10)}M_{(10,11)}M_{(7,8)}}.$$

Let us also write $B$-invariants of the first series,
$\Psi_1=\big\{(5,9),(8,12)\big\}$,
$$A_{(5,9)}=\frac{L_{(5,9)}L_{(6,8)}}{L_{(5,8)}M_{(3,6)}M_{(6,9)}},\quad
A_{(8,12)}=\frac{L_{(8,12)}L_{(9,11)}}{L_{(8,11)}M_{(5,9)}M_{(9,12)}}.$$
}

\thenv{Lemma\label{B_invariant}}{\emph{Rational functions}
$B_{\psi}$, $\psi\in\Psi_2$, \emph{are invariant under the adjoint
action of} $B$.}

\textsc{Proof.} Similarly to the case $A_{\psi}$, it is sufficient
to prove that $B_{\psi}$ are $T$-invariant, where $T$ is the set of
nondegenerate diagonal matrices. We prove the statment for the case
$r_s=r_t$. The others cases are similar.

If $\psi\in\Psi_2$ and $r_s=r_t$, then there exist roots
$\xi_1\in\Phi$ and $\gamma_1,\ldots,\gamma_5\in S$ describing
before. We have
$$L_{\psi}=D_{I_1,\emptyset}^{\{R_t+1\},J'}\!\!,\quad
L_{\xi_1}=D_{\{R_{s-1}\},I'}^{J_1,\emptyset},\quad
M_{\gamma_1}=\mathbb{X}_{\{R_{s-1}\}}^{\{R_{s-1}+1\}},\quad
M_{\gamma_3}=\mathbb{X}_{\{R_t\}}^{\{R_t+1\}},$$
$$M_{\gamma_5}=\mathbb{X}_{\{R_{s-1}+1,R_{s-1}+2,\ldots,R_{t-1}\}}^%
{\{R_s+1,R_s+2,\ldots,R_t\}},\quad \prod_{\mu\prec\gamma_4}M_{\mu}=
\mathbb{X}_{I_1}^{J_1},$$
$$\prod_{\mu'\prec\mu\prec\gamma_4}M_{\mu}=
\mathbb{X}_{I'}^{J'},\quad\prod_{\mu\prec\gamma_5}M_{\mu}=
\mathbb{X}_{\{R_{s-1}+2,R_{s-1}+3,\ldots,R_{t-1}\}}^%
{\{R_s+1,R_s+2,\ldots,R_t-1\}},$$ where
$$I_1=\big\{\mathrm{row}(\gamma_4)+1,\mathrm{row}(\gamma_4)+2,
\ldots,R_{t-1}\big\},\
J_1=\big\{R_s+1,R_s+2,\ldots,\mathrm{col}(\gamma_4)-1\big\},$$
$$I'=I_1\setminus\{\mathrm{row}(\gamma_4)+1\},\mbox{ and }
J'=J_1\setminus\{\mathrm{col}(\gamma_4)-1\}.$$ Note that the number
of elements in the sets $I_1$ and $J_1$ are the same because
$\mathbb{X}_{I_1\cup\{\mathrm{row}(\gamma_4)\}}%
^{J_1\cup\{\mathrm{col}(\gamma_4)\}}=M_{\gamma_4}$ is a square
minor.

\medskip
Then
$$\mathrm{Ad}_t\,B_{\psi}=\mathrm{Ad}_t\,
\frac{\displaystyle L_{\xi_1}L_{\psi}\cdot
\prod_{\mu\prec\gamma_5}M_{\mu}}{\displaystyle
D\cdot\prod_{\mu\prec\gamma_4}M_{\mu}\cdot
\prod_{\mu'\prec\mu\prec\gamma_4}M_{\mu'}}=\frac{\displaystyle\frac{t_{R_{s-1}}t_{I'}}{
t_{J_1}}L_{\xi_1}\cdot\frac{t_{I_1}}{t_{R_t+1}t_{J'}}L_{\psi}
}{\displaystyle\frac{t_{R_{s-1}}}{
t_{R_{s-1}+1}}M_{\gamma_1}\cdot\frac{t_{R_t}}{t_{R_t+1}}M_{\gamma_3}}
\times$$
$$\times
\frac{\displaystyle\frac{t_{R_{s-1}+2}t_{R_{s-1}+3}\ldots
t_{R_{t-1}}}{t_{R_s+1}t_{R_s+2}\ldots t_{R_t-1}}
\prod_{\mu\prec\gamma_5}M_{\mu}}{\displaystyle\frac{t_{R_{s-1}+1},t_{R_{s-1}+2},
\ldots,t_{R_{t-1}}}{t_{R_s+1},
t_{R_s+2},\ldots,t_{R_t}}M_{\gamma_5}\cdot
\frac{t_{I_1}}{t_{J_1}}\prod_{\mu\prec\gamma_4}M_{\mu}\cdot
\frac{t_{I'}}{t_{J'}}\prod_{\mu'\prec\mu\prec\gamma_4}M_{\mu}}=B_{\psi}.$$
So $B_{\psi}$ is $T$-invariant.~$\Box$

\medskip
Denote
$$\mathcal{X}=\sum_{\psi\in\Psi}c_{\psi}E_{\psi}+
\sum_{\xi\in(S\cup\Phi)\setminus\Psi}E_{\xi}.$$ We will show below
that a $B$-orbit in general position has a representative in
$\mathcal{X}$.

Let $\pi:K(\mathfrak{m})^B\rightarrow K(\mathcal{X})$ be the
restriction map
\begin{equation}
\pi(x_{\xi})=x_{\xi}|_{\mathcal{X}}=\left\{\begin{array}{ll}
c_{\xi}&\mbox{if }\xi\in\Psi,\\
1&\mbox{if }\xi\in(S\cup\Phi)\setminus\Psi,\\
0&\mbox{if }\xi\not\in S\cup\Phi,
\end{array}\right.\label{pi}
\end{equation}
where the field $K(\mathcal{X})$ is a field of fractions for the
polynomial algebra $K[\mathcal{X}]$ of variables $c_{\psi}$,
$\psi\in\Psi$.

\thenv{Proposition\label{A-B_independ}}{\emph{The system of rational
functions}
$$\big\{A_{\psi_1},\ B_{\psi_2}:\ \psi_1\in\Psi_1,\
\psi_2\in\Psi_2\big\}$$
\emph{is algebraically independent}.}

\textsc{Proof.} It is sufficient to prove that the system
$$\big\{A_{\psi_1}|_{\mathcal{X}},\ B_{\psi_2}|_{\mathcal{X}},\
\psi_1\in\Psi_1,\ \psi_2\in\Psi_2\big\}$$ is algebraically
independent.

Let us number roots $\psi\in\Psi$ from the bottom up in columns from
left to right.

\thenv{Example}{Roots $\psi\in\Psi$ are numbered, roots from $S$ and
from $\Phi\setminus\Psi$ are labeled by the symbols $\otimes$ and
$\times$ as before.
\begin{center}\refstepcounter{diagram}
{\begin{tabular}{|p{0.1cm}|p{0.1cm}|p{0.1cm}|p{0.1cm}|p{0.1cm}|p{0.1cm}
|p{0.1cm}|p{0.1cm}|p{0.1cm}|p{0.1cm}|p{0.1cm}|p{0.1cm}|c}
\multicolumn{2}{l}{{\small 1\quad 2}}&\multicolumn{2}{l}{{\small
3\quad 4\!}}&\multicolumn{2}{l}{{\small 5\quad
6\!}}&\multicolumn{2}{l}{{\small 7\quad
8\!}}&\multicolumn{2}{l}{{\small
9\ \  10\!\!}}&\multicolumn{2}{l}{{\small 11\  12\!\!\!}}\\
\cline{1-12} \multicolumn{3}{|l|}{1}&&&$\otimes$&&&&&&&{\small 1}\\
\cline{4-12} \multicolumn{3}{|c|}{1}&&$\otimes$&&&&&&&&{\small 2}\\
\cline{4-12} \multicolumn{3}{|r|}{\qquad\ \ 1}&$\otimes$&&&&&&&&&{\small 3}\\
\cline{1-12} \multicolumn{3}{|c|}{}&\multicolumn{4}{|l|}{1}&$\times$&2&3&&&{\small 4}\\
\cline{8-12} \multicolumn{3}{|c|}{}&\multicolumn{4}{|l|}{\quad\ 1}&$\times$&1&$\otimes$&&&{\small 5}\\
\cline{8-12} \multicolumn{3}{|c|}{}&\multicolumn{4}{|l|}{\qquad\ \ 1}&$\times$&$\otimes$&&&&{\small 6}\\
\cline{8-12} \multicolumn{3}{|c|}{}&\multicolumn{4}{|r|}{1}&$\otimes$&&&&&{\small 7}\\
\cline{4-12} \multicolumn{7}{|c|}{}&\multicolumn{3}{|l|}{1}&5&6&{\small 8}\\
\cline{11-12} \multicolumn{7}{|c|}{}&\multicolumn{3}{|c|}{1}&4&$\otimes$&{\small 9}\\
\cline{11-12} \multicolumn{7}{|c|}{}&\multicolumn{3}{|r|}{1}&$\otimes$&&{\small 10}\\
\cline{8-12} \multicolumn{10}{|c|}{}&\multicolumn{2}{|l|}{1}&{\small 11}\\
\multicolumn{10}{|c|}{}&\multicolumn{2}{|r|}{1}&{\small 12}\\
\cline{1-12} \multicolumn{12}{c}{Diagram \arabic{diagram}}\\
\end{tabular}}
\end{center}}

We prove the lemma by induction on number of roots in $\Psi$. Since
for the first root $\psi\in\Psi$ we have $A_{\psi}=c_{\psi}$ if
$\psi\in\Psi_1$ and $B_{\psi}=c_{\psi}$ if $\psi\in\Psi_2$, the
system consisting of one $B$-invariant is algebraically independent.
Therefore the base of induction is obvious. Suppose that the system
$\big\{A_{\psi'},B_{\psi''}\big\}$, where $\psi'$ and $\psi''$ are
all roots in $\Psi$ such that numbers of $\psi'$ and $\psi''$ are
less than $k$, is algebraically independent. Let $\psi\in\Psi$ be a
root with the number $k$. For any $\xi\in S$ and $\varphi\in\Phi$ we
have
$$M_{\xi}|_{\mathcal{X}}=1,\quad
L_{\varphi}|_{\mathcal{X}}=\left\{\begin{array}{ll}
c_{\varphi}&\mbox{if }\varphi\in\Psi,\\
1&\mbox{if }\varphi\in\Phi\setminus\Psi.\\
\end{array}\right.\ $$
Denote
$$\widetilde{c}_{\varphi}=\left\{\begin{array}{ll}
c_{\varphi}&\mbox{if }\varphi\in\Psi,\\
1&\mbox{if }\varphi\in\Phi\setminus\Psi\\
\end{array}\right.$$
for any $\varphi\in\Phi$. Then
\begin{equation}
A_{\psi}|_{\mathcal{X}}=\left\{\begin{array}{ll}
\displaystyle\frac{c_{\psi}\widetilde{c}_{\xi_2}}%
{\widetilde{c}_{\xi_1}\widetilde{c}_{\xi_3}}&\mbox{if
}\xi_3\in\Phi,\\
\\
\displaystyle\frac{c_{\psi}\widetilde{c}_{\xi_2}}%
{\widetilde{c}_{\xi_1}}&\mbox{if
}\xi_3\in S;\\
\end{array}\right.\qquad
B_{\psi}|_{\mathcal{X}}=\left\{\begin{array}{ll} \displaystyle
c_{\psi}\widetilde{c}_{\xi_1}&\mbox{if }r_s=r_t,\\
\\
\displaystyle\frac{c_{\psi}\widetilde{c}_{\xi_1}}{\widetilde{c}_{\xi_2}}&
\mbox{if }r_s<r_t,\\
\\
\displaystyle\frac{c_{\psi}\widetilde{c}_{\xi_1}}{\widetilde{c}_{\xi_3}}&
\mbox{if }r_s>r_t,\\
\end{array}\right.\label{A_psi|_X}
\end{equation}
where $\xi_1,\xi_2,\xi_3$ are the mentioned above roots
corresponding to $\psi$. Note that the roots $\xi_1,\xi_2,\xi_3$ are
below of the root $\psi$ or the numbers of columns for
$\xi_1,\xi_2,\xi_3$ are less than $\mathrm{col}(\psi)$ for the both
cases $\psi\in\Psi_1$ and $\psi\in\Psi_2$. Therefore the numbers of
$\xi_1,\xi_2,\xi_3$ are less than $k$. Hence there is no the
variable $c_{\psi}$ in any function of the system
$\big\{A_{\psi'},B_{\psi''}\big\}$, where numbers of $\psi'$ and
$\psi''$ are less than $k$. It means that the image~(\ref{A_psi|_X})
of the $B$-invariant corresponding to $\psi$ does not depend on
$\big\{A_{\psi'},B_{\psi''}\big\}$, where numbers of $\psi'$ and
$\psi''$ are less than $k$. So we get the statement of the
proposition.~$\Box$


\sect{The canonical representative for $B$-orbits in general
position}

In this section we show that a $B$-orbit in general position has a
unique repre\-sen\-tative in the set
$$\mathcal{X}=\sum_{\psi\in\Psi}c_{\psi}E_{\psi}+
\sum_{\xi\in(S\cup\Phi)\setminus\Psi}E_{\xi},\quad c_{\psi}\neq0.$$

\thenv{Theorem\label{Exist_representative_for_B}}{\emph{There exists
an nonempty open set $U\subset\mathfrak{m}$ such that for any $x\in
U$ there is $g\in B$ satisfying} $\mathrm{Ad}_gx\in\mathcal{X}$.}

\textsc{Proof.} Let $T$ be the set of nondegenerate diagonal
matrices. By Theo\-rem~\ref{Exist_of_representative} and using
$B=T\ltimes N$ it is sufficient to show that $W=U$ and there is
$t\in T$ such that for any $x\in\mathcal{Y}\cap W$ the $T$-orbit of
$x$ intersects $\mathcal{X}$ at a unique point, where $W$ is
nonempty open from Theorem~\ref{Exist_of_representative}.

We describe an algorithm of bringing $x\in\mathcal{Y}\cap W$ to the
canonical form. Denote
$$h_i(b)=\mathrm{diag}(\underbrace{1,\ldots,1}_{i-1\mbox{ \small
units}},b,\underbrace{1,\ldots,1}_{n-i\mbox{ \small units}}).$$ Then
the adjoin action of $h_i(b)$ on $y\in\mathfrak{m}$ consists of two
transformations:
\begin{itemize}
\item[1)] the row $i$ is multiplied by $b$,
\item[2)] the column $i$ is multiplied by $b^{-1}$.
\end{itemize}
Let $a_1=1$ and $a_1<a_2<\ldots<a_p$ are numbers of blocks in
$\mathfrak{r}$ such that $1<r_{a_2}<\ldots<r_{a_p}$ and $p$ is
maximal.

Let the map $\omega_{i,j}:\mathfrak{m}\rightarrow K$ take each
$y\in\mathfrak{m}$ to the element in the intersection of the $i$th
row and the $j$th column of $y$. Denote $\omega_{\xi}=\omega_{i,j}$
if $\xi=(i,j)$.

We need making units in positions of the matrix $x$ corresponding to
$\xi$ for $\xi\in(S\cup\Phi)\setminus\Psi$. The algorithm consists
of 5 steps.

\medskip
\textsc{Step 1.} Consider the first column of the block
$X_{a_2,a_2+1}$, it is the column $R_{a_2}+1$. By the definition of
$a_1,a_2,\ldots,a_p$ it follows that $r_2=r_3=\ldots=r_{a_2-1}=1$.
Therefore there are not roots of $\Phi$ in the rows
$1,2,\ldots,R_{a_2-1}$. The pair
$$q=\big((R_{a_2-1},R_{a_2-1}+1),(R_{a_2},R_{a_2}+1)\big)$$ is
admissible. Hence the root $\xi=(R_{a_2-1}+1,R_{a_2}+1)\in\Phi$
correspond to $q$ and the row $R_{a_2-1}+1$ is the first row with a
root of $\Phi$. So $\xi\not\in\Psi$. Similarly, there are not roots
of $\Psi$ in the column $R_{a_2}+1$.

Note that since $x\in\mathcal{Y}$, we have $\omega_{\xi}(x)\neq0$
for any $\xi$ of the extended base. Therefore for any $b\in T$ and
$\xi\in S\cup\Phi$ we have $\omega_{\xi}(\mathrm{Ad}_b\,x)\neq0$.

The adjoint action of the element
$$g_1=\prod_{i=R_{a_2-1}+1}^{R_{a_2}}h_{i}(\omega_{i,
R_{a_2}+1}(x)^{-1})$$ on $x$ makes units in the positions
$(R_{a_2-1}+i,R_{a_2}+1)$ for $i=1,2,\ldots,r_{a_2}$, i.e. the
action of $g_1$ makes units in the first column of $X_{a_2,a_2+1}$.
The numbers $1$ label the cells in the diagram
Example~\ref{Ex-213142-ex}, where we make units by $g_1$.

\thenv{Example\label{Ex-213142-ex}}{Let the diagonal sizes are
$(2,1,3,1,4,2)$, then $p=3$ and $a_1=1$, $a_2=3$, $a_3=5$. Note that
if $\xi\in S\cup\Phi$ belongs to $\Psi$, then $\xi$ is labeled by
the symbol $\boxtimes$. If $\xi\in(S\cup\Phi)\setminus\Psi$, then
$\xi$ is labeled by numbers of steps.
\begin{center}\refstepcounter{diagram}
{\begin{tabular}{|p{0.1cm}|p{0.1cm}|p{0.1cm}|p{0.1cm}|p{0.1cm}|p{0.1cm}
|p{0.1cm}|p{0.1cm}|p{0.1cm}|p{0.1cm}|p{0.1cm}|p{0.1cm}|p{0.1cm}|c}
\multicolumn{2}{l}{{\small 1\quad 2}}&\multicolumn{2}{l}{{\small
3\quad 4\!}}&\multicolumn{2}{l}{{\small 5\quad
6\!}}&\multicolumn{2}{l}{{\small 7\quad
8\!}}&\multicolumn{2}{l}{{\small
9\ \  10\!\!\!}}&\multicolumn{3}{l}{{\small 11\  12\ 13\!\!\!}}\\
\cline{1-13} \multicolumn{2}{|l|}{1}&&&$\!2$&&&&&&&&&{\small 1}\\
\cline{3-13} \multicolumn{2}{|r|}{1}&$\!2$&&&&&&&&&&&{\small 2}\\
\cline{1-13} \multicolumn{2}{|r|}{}&1&$\!2$&&&&&&&&&&{\small 3}\\
\cline{3-13} \multicolumn{3}{|c|}{}&\multicolumn{3}{|l|}{1}&$\!1$&&$\!\boxtimes$&$\!3$&&&&{\small 4}\\
\cline{7-13} \multicolumn{3}{|c|}{}&\multicolumn{3}{|c|}{1}&$\!1$&&$\!3$&&&&&{\small 5}\\
\cline{7-13} \multicolumn{3}{|c|}{}&\multicolumn{3}{|r|}{1}&$\!1$&&&&&&&{\small 6}\\
\cline{4-13} \multicolumn{6}{|c|}{}&1&$\!3$&&&&&&{\small 7}\\
\cline{7-13} \multicolumn{7}{|c|}{}&\multicolumn{4}{|l|}{1}&$\!\boxtimes$&$\!\boxtimes$&{\small 8}\\
\cline{12-13} \multicolumn{7}{|c|}{}&\multicolumn{4}{|l|}{\,\quad1}&$\!\boxtimes$&$\!\boxtimes$&{\small 9}\\
\cline{12-13} \multicolumn{7}{|c|}{}&\multicolumn{4}{|c|}{\quad1}&$\!4$&$\!5$&{\small 10}\\
\cline{12-13} \multicolumn{7}{|c|}{}&\multicolumn{4}{|r|}{1}&$\!4$&&{\small 11}\\
\cline{8-13} \multicolumn{11}{|c|}{}&\multicolumn{2}{|l|}{1}&{\small 12}\\
\multicolumn{11}{|c|}{}&\multicolumn{2}{|r|}{1}&{\small 13}\\
\cline{1-13} \multicolumn{12}{c}{Diagram \arabic{diagram}}\\
\end{tabular}\label{Ex-213242}}
\end{center}}

\textsc{Step 2.} Any row $1,2,\ldots,R_{a_2-1}$ of a diagram
contains a single symbol $\otimes$. The adjoint action of
$$g_2=\prod_{i=1}^{R_{a_2-1}}h_{i}(b_i)$$
with a suitable choice of $b_i$ on $\mathrm{Ad}_{g_1}x$ makes units
in the positions corresponding to roots of $S$ in rows
$1,2,\ldots,R_{a_2-1}$ and $g_2$ leaves the rest part of
$\mathrm{Ad}_{g_1}x$ unchanged. We have
$b_i=\omega_{\xi_i}(\mathrm{Ad}_{g_1}x)^{-1}$ for the root $\xi_i$
in the $i$th row. The step 2 makes units in the positions that are
label by 2 in Diagram~\ref{Ex-213242}. Now we have that the first
$R_2+1$ columns have the canonical form.

Take $k:=3$.

\medskip
\textsc{Step 3} provides the canonical form in the columns
$R_{a_{k-1}}+2,\ldots,R_{a_k}$.

Let $j$ be any column of $R_{a_{k-1}}+2,\ldots,R_{a_k}$. Let us show
that if there exists $\psi\in\Phi$ in this column, then
$\psi\in\Psi$. There are two cases, the first case is $j=R_t+1$ for
some $t$, where $a_{k-1}+1<t\leqslant a_k$, and the second one is
$j\neq R_t+1$ for any $t$.

For the first case $j=R_t+1$ by the definition of $a_1,\ldots,a_p$
we have $r_{t}\leqslant r_{a_{k-1}}$. By Lemmas~\ref{Lemma1}
and~\ref{Lemma2} there exists a root of the base in every row
$R_{a_{k-1}-1}+1,R_{a_{k-1}-1}+2,\ldots,R_{a_{k-1}}$. Therefore the
root $(R_{a_{k-1}-1},R_{a_{k-1}-1}+1)$ and every of these roots
except the root in the row $R_{a_{k-1}-1}+1$ make an admissible
pair. Hence there are $r_{a_{k-1}}-1$ roots of $\Phi$ in the row
$R_{a_{k-1}-1}+1$. Since $r_{t}\leqslant r_{a_{k-1}}$ there are no
more that $r_t-1$ roots of $\Phi$ in the column $j=R_t+1$. Therefore
by the definition of $\Psi_2$ every root of $\Phi$ in the column $j$
belongs to $\Psi_2$.

Consider the second case. If $\psi=(i,j)$ is in $\Phi$ for some $i$
and $R_t+1<j\leqslant R_{t+1}$ for some $t$, then let us prove that
$\psi\in\Psi_1$. It is sufficient to prove that
$\xi_1=(i,R_t+1)\in\Phi$, $\xi_2=(i+1,R_t+1)\in\Phi$, and
$\xi_3=(i+1,j)\in S\cup\Phi$. Since $(i,j)\in\Phi$, then there is
the first root $\gamma_1\in S$ of the admissible pair for $\psi$ in
the column $i$. Then $\big(\gamma_1,(R_t,R_t+1)\big)$ is an
admissible pair and the corresponding root $(i,R_t+1)$ is in $\Phi$.
Since $r_t\leqslant r_{k-1}$, by Lemmas~\ref{Lemma1}
and~\ref{Lemma2} we have that there is a root of $S$ in every column
$R_{t-1}+1,R_{t-1}+2,\ldots,R_t$, i.e. in every column above the
block $r_t$ in $\mathfrak{r}$. Therefore there exists $\gamma_2\in
S$ in the column $i+1$, then $\big(\gamma_2,(R_t,R_t+1)\big)$ is
admissible and $\xi_2\in\Phi$. Similar reasoning shows $\xi_3\in
S\cup\Phi$. Thus we have $\psi\in\Psi$.

So we need constructing units in columns
$R_{a_{k-1}}+2,\ldots,R_{a_k}$ only in the positions $\varphi$,
where $\varphi\in S$. The action of
$$g_{2k-3}=\prod_{j=R_{a_{k-1}}+2}^{R_{a_k}}h_{j}(b_j),\quad
b_j=\omega_{\varphi_j}(\mathrm{Ad}_{g_{2k-4}\ldots g_1}x)$$ on
$\mathrm{Ad}_{g_{2k-4}\ldots g_1}x$, where $\varphi_j\in S$ is the
root in the $j$th column, constructs units in positions
corresponding to roots in $S$ in columns $j$,
$R_{a_{k-1}}+1<j\leqslant R_{a_k}$. If there is no a root of $S$ in
the $j$th column, then $h_j$ is equal to the identity matrix. The
action of $g_{2k-3}g_{2k-4}\ldots g_1$ multiplies the rows
$R_{a_{k-1}}+2,\ldots,R_{a_k}$ and the same columns on $b_i$ and
therefore it makes the canonical view in the columns
$1,2,\ldots,R_{a_k}$.

In Example~\ref{Ex-213142-ex} Step 3 constructs units in positions
that was labeled by 3.

\medskip
\textsc{Step 4} is performed if $k\neq p$ or $R_{a_p}<n$.

Consider all roots of $(S\cup\Phi)\setminus\Psi$ in the column
$R_k+1$ to the right of the $r_k$th block in $\mathfrak{r}$. We show
in the previous step that for any $m$, $1<m<p$, there are exactly
$a_m-1$ roots of $\Phi$ in the row $R_{a_{m-1}}+1$. Therefore the
maximal number of roots of $\Phi$ in a row $i$, where $i\leqslant
R_{a_{k-1}}$, is $r_{a_{k-1}}+1$. Hence if $1\leqslant
i<r_{a_{k-1}}$, then $\xi=(R_{a_k-1}+i,R_{a_k}+1)\in\Psi$ and if
$r_{a_{k-1}}\leqslant i\leqslant r_{a_k}$, then
$\xi=(R_{a_k-1}+i,R_{a_k}+1)\in\Phi\setminus\Psi$. By
Lemmas~\ref{Lemma1} and~\ref{Lemma2} there exist $r_{a_{k-1}}$ roots
of the base in the first $r_{a_{k-1}}$ columns above the $r_k$th
block in $\mathfrak{r}$. Besides, the root of $S$ in the last column
$R_{a_k-1}+r_{a_{k-1}}$ of these columns is in the row
$R_{a_{k-1}-1}+1$. Therefore there is not a root of $\Phi$ in the
column $R_{a_k-1}+r_{a_{k-1}}$. Besides there are not roots of the
extended base in columns $R_{a_k-1}+r_{a_{k-1}}+1,\ldots,R_{a_k}$.

The action of element
$$g_{2k-2}=\left(
\prod_{i=R_{a_k-1}+r_{a_{k-1}}+1}^{R_{a_k}} h_{i}(b_i)\right)\cdot
h_{R_{a_k}+1}(b_{R_{a_k}+1})$$
$$\mbox{for }b_{R_{a_k}+1}=\omega_{R_{a_k-1}+r_{a_{k-1}},R_{a_k}+1}
(\mathrm{Ad}_{g_{2k-1}\ldots g_1}x)$$
$$\mbox{and }b_{i}=\omega_{i,R_{a_k}+1}(\mathrm{Ad}_{g_{2k-1}\ldots
g_1}x)^{-1}$$ makes units in the positions
$(R_{a_k-1}+i,R_{a_k}+1)$, where $r_{a_{k-1}}\leqslant i\leqslant
R_{a_k}$ and does not change the rows $1,2,\ldots,R_{a_k-1}$ of the
matrix $\mathrm{Ad}_{g_{2k-1}\ldots g_1}x$.

We increase the variable $k$ by 1 and go to the Step 3 until $k=p$.

\medskip
\textsc{Step 5} is performed if $R_{a_p}+1<n$. Using the similar
reasoning, one has that any root of $\Phi$ in columns $j$,
$j>R_{a_p}+1$, belongs to $\Psi$.

We make units in positions corresponding to roots of $S$ in the
columns $R_{a_p}+2,\ldots,n$ as in Step 3.~$\Box$

\thenv{Theorem}{\emph{The restriction map}~(\ref{pi}) \emph{is a
bijection from $K(\mathfrak{m})^B$ to} $K(\mathcal{X})$.}

\textsc{Proof.} Let us show that $\pi$ is an injection. Suppose
$\pi(f)=0$ for some $f\in K(\mathfrak{m})^B$. It means
$f(\mathcal{X})=0$. By Theorem~\ref{Exist_representative_for_B} we
have $\overline{\mathrm{Ad}_B\mathcal{X}}=\mathfrak{m}$. Therefore
since $f$ is a $B$-invariant, we obtain
$$f(\mathrm{Ad}_B\mathcal{X})=f(\mathcal{X})=0.$$
Hence $f\equiv0$ and $\mathrm{Ker}(\pi)=\{0\}$.

Let us check that any variable $c_{\psi}$, $\psi\in\Psi$, has a
preimage in $K(\mathfrak{m})^B$ under the action of $\pi$. We use
the numbering $\psi_1,\psi_2,\ldots,\psi_k$, $k=|\Psi|$, for the
roots of $\Psi$ from the proof of Theorem~\ref{A-B_independ}. Let us
prove the theorem by induction on the number of roots from $\Psi$.

The base of induction is showed as follows. If $\psi_1\in\Psi_1$,
then by definition of $A_{\psi_1}$ there exist roots
$\xi_1,\xi_2,\xi_3\in S\cup\Phi$. Since the number of $\psi_1$ is 1
and $\xi_1,\xi_2,\xi_3$ are below or to the left of $\psi$, then
$\xi_1,\xi_2,\xi_3$ do not belong to $\Psi$. Therefore,
$c_{\psi_1}=\pi(A_{\psi_1})$. Similarly, if $\psi_1\in\Psi_2$, then
$c_{\psi_1}=\pi(B_{\psi_1})$.

Suppose that the statement is true for any number less than $m$. If
$\psi_m\in\Psi_1$, then
$$A_{\psi_m}|_{\mathcal{X}}=
\frac{c_{\psi_m}c_{\xi_2}}{c_{\xi_1}c_{\xi_3}}\mbox{ if }
\xi_3\in\Phi\mbox{ or }
A_{\psi_m}|_{\mathcal{X}}=\frac{c_{\psi_m}c_{\xi_2}}{c_{\xi_1}}
\mbox{ if }\xi_3\in S$$ for the roots $\xi_1,\xi_2,\xi_3$ defined
for $\psi_m$. Since $\xi_1,\xi_2,\xi_3$ are under and to the left of
$\psi_m$, then the numbers of $\xi_1,\xi_2,\xi_3$ are less than the
number of $\psi_m$. Therefore by the inductive assumption there
exist $B$-invariant rational functions $f,g,h$ such that
$$c_{\xi_1}=\pi\big(f(A_{\phi},B_{\varphi})_{\phi,\varphi\in\Psi}\big),$$
$$c_{\xi_2}=\pi\big(g(A_{\phi},B_{\varphi})_{\phi,\varphi\in\Psi}\big),$$
$$c_{\xi_3}=\pi\big(h(A_{\phi},B_{\varphi})_{\phi,\varphi\in\Psi}\big).$$
Then in the case $\xi_3\in\Phi$ we have
$$c_{\psi_m}=\pi\left(\frac{A_{\psi_m}f(A_{\phi},B_{\varphi})
h(A_{\phi},B_{\varphi})}{g(A_{\phi},B_{\varphi})}\right).$$ If
$\xi_3\in S$, then
$$c_{\psi_m}=\pi\left(\frac{A_{\psi_m}f(A_{\phi},B_{\varphi})
}{g(A_{\phi},B_{\varphi})}\right).$$ So $c_{\psi_m}$ has a preimage.

The case $\psi_m\in\Psi_2$ is proved similarly.~$\Box$

\thenv{Corollary\label{field_B-invariant}}{\emph{The field of
invariants $K(\mathfrak{m})^B$ is generated by rational functions
$A_{\psi_1}$ and} $B_{\psi_2}$, $\psi_1\in\Psi_1$,
$\psi_2\in\Psi_2$.}

Using Corrolary of Theorem 2.3 of~\cite{PV} one has

\thenv{Corollary}{\emph{The dimension of a $B$-orbit in general
position equals $$\dim\,\mathfrak{m}-|\Psi|.$$ The rational
functions $A_{\psi_1}$ and} $B_{\psi_2}$, $\psi_1\in\Psi_1$,
$\psi_2\in\Psi_2$, \emph{separate $B$-orbits in general position}.}

\medskip
\textbf{Acknowledgments.}  \emph{I am very grateful to Anna Melnikov
and Alexander Panov for discussions.}

\textsc{Department of Mathematics, University of Haifa, Israil}\\
\textsc{Department of Mechanics and Mathematics, Samara University,
Russia}\\ \emph{E-mail address}: \verb"berlua@mail.ru"


\begin{thebibliography}{99}
\bibitem[GG]{GG}
М. Goto and F. Grosshans, Semisimple Lie algebras, Lect. Notes in
Pure Appl. Math., vol. 38 (1978).
\bibitem[K]{K}
H. Kraft, Geometrische Methoden in der Invariantentheorie, Friedr.
Vieweg and Sohn, Braunschweig/Wiesbaden (1985).
\bibitem[PS]{PS}
A. N. Panov and V. V. Sevostyanova, Regular $N$-orbits in the
nilradical of a parabolic subalgebra, \emph{Vestnik SamGU},
\textbf{7}(57) (2007), pp. 152--161. See also
http://arxiv.org/abs/1203.2754.
\bibitem[PV]{PV}
V. L. Popov and E. B. Vinberg, Invariant theory, in:
\emph{Progressin Science and Technology}, VINITI, Moscow (1989), pp.
137--309.
\bibitem[R]{R}
R. W. Richardson, Conjugacy classes in parabolic subgroups of
semisimple  algebraic groups, \emph{Bull. London Math. Soc.}
\textbf{6} (1974), pp. \mbox{21--24.}
\bibitem[S1]{S1}
V. V. Sevostaynova, The invariant field of the adjoint action of the
unitriangular group in the nilradical of a parabolic subalgebra.
\emph{Zapiski nauchn. seminarov POMI}, vol. 375, 2010, pp. 167--194.
(English translation: \emph{Journal of Math. Sciences}, vol. 171,
No. 3, 2010, pp. 400--415.) See also http://arxiv.org/abs/1203.3000.
\bibitem[S2]{S2}
V. V. Sevostyanova, The invariant algebra of the adjoint action of
the unitriangular group in the nilradical of a parabolic subalgebra.
Vest. SamGU \textbf{2}(76) (2010), pp. 72--83. See also
http://arxiv.org/abs/1203.4899.
\bibitem[S3]{S3}
V. V. Sevostyanova, The algebra of invariants for the adjoint action
of the unitriangular group. http://arxiv.org/abs/1605.00800

\end{thebibliography}
\end{document}